# Periodicity of Hermitian K-theory and Milnor's K-groups

## by Max Karoubi

**Introduction**

One purpose of this paper is to improve the periodicity theorem of Hermitian K-theory proved in [6], using some recent results by A..J. Berrick and the author [1]. For instance, if $1/2 \in A$ we show that the "higher Witt groups" $_\varepsilon W_n(A)$, defined for a ring A with antiinvolution ($\varepsilon$ being the sign of symmetry), are periodic of period $\alpha \equiv 0$ mod. 4, with a control of the 2-torsion depending on $\alpha$. More precisely, we have homomorphisms

$$\beta : {}_\varepsilon W_n(A) \longrightarrow {}_\varepsilon W_{n+\alpha}(A) \text{ and } \beta' : {}_\varepsilon W_{n+\alpha}(A) \longrightarrow {}_\varepsilon W_n(A)$$

such that $\beta.\beta'$ and $\beta'.\beta$ are multiplication by $2.8^r$, if $\alpha = 8r$ or $\alpha = 8r - 4$ (with $r > 0$).
If $1/2 \notin A$, we also show the existence of analogous homomorphisms whose composites are multiplication by 64 if $\alpha = 4$ and by 512 if $\alpha = 8$ (a result which is probably not optimal).

If A is *any* commutative ring and if $1/2 \in A$, this periodicity phenomenon is related to a new decreasing filtration $(F_n)$ of the classical Witt ring $W(A)$ [denoted by $_1W_0(A)$ with our general conventions]. More precisely, $F_{2r}$ is the image of the "periodicity homomorphism":

$$_{(-1)^r}W_{2r}(A) \longrightarrow W(A).$$

We write $F_1 = I$ for the augmentation ideal, kernel of the "rank map"

$$W(A) \longrightarrow k_0(A) = H^0(\mathbf{Z}/2, K_0(A)).$$

There is a more general rank map[1]

$$_{(-1)^r}W_{2r}(A) \longrightarrow k_{2r}(A) = H^0(\mathbf{Z}/2, K_{2r}(A))$$

whose kernel is denoted by $_{(-1)^r}W_{2r}(A)'$. Finally, we define $F_{2r+1}$ as the image of the composite

---
[1] where $K_n(A)$ denotes the Quillen K-groups in general.

$$_{(-1)^r}W_{2r}(A)' \longrightarrow {}_{(-1)^r}W_{2r}(A) \longrightarrow W(A)$$

where the first map is an inclusion. With these definitions, one can prove (at least for regular Noetherian rings) that $F_n F_p \subset F_{n+p}$ and therefore $F_n \subset I^n$ (a strict inclusion in general).

If A is a field of characteristic $\neq 2$ with trivial involution, the group $k_0(A)$ is reduced to $\mathbf{Z}/2$. This periodicity of the higher Witt groups $W_*(A)$ and the filtration above are then related to Milnor's K-groups mod. 2 denoted by $k_{2r}^M(A)$, thanks to the work of Voevodsky [7]. For instance, if $k_{2r}^M(A)$ is not trivial, we have $_{(-1)^s}W_{2s}(A) \neq 0$ for all $s \leq r$ (more general results are stated in the second section of this paper).

## 1. Periodicity of Hermitian K-theory.

**1.1.** We start by recalling briefly some notations and results which are well known by the experts. Let A be a ring with an antiinvolution $a \mapsto \bar{a}$ and let $\varepsilon$ be an element of the center of A such that $\varepsilon \bar{\varepsilon} = 1$. Associated to those data, there is a space $_\varepsilon\mathcal{K}(A)$ whose homotopy groups are the "higher Hermitian K-groups", denoted by $_\varepsilon L_n(A)$.

More precisely, $_\varepsilon L_0(A)$ is the Witt-Grothendieck group of the category of finitely generated projective modules provided with non degenerate $\varepsilon$-Hermitian forms. The groups $_\varepsilon L_n(A)$ for $n > 0$ are also homotopy groups $\pi_n(B_\varepsilon O(A)^+)$, as explained below.

First, we define $_\varepsilon O_{n,n}(A)$ as the group whose elements are given by 2n x 2n matrices

$$M = \begin{pmatrix} a & b \\ c & d \end{pmatrix}$$

such that $M^*M = MM^* = 1$, the identity matrix, where

$$M^* = \begin{pmatrix} {}^t\bar{d} & \varepsilon {}^t\bar{b} \\ \bar{\varepsilon}\,{}^t\bar{c} & {}^t\bar{a} \end{pmatrix}$$

The infinite $\varepsilon$–orthogonal group

$$_\varepsilon O(A) = \lim {}_\varepsilon O_{n,n}(A)$$

has a commutator subgroup which is perfect. We can perform Quillen's + construction on its classifying space in order to kill the commutator subgroup. This enables us to define the higher homotopy groups $\pi_n(B_\varepsilon O(A)^+)$ introduced above. The space $_\varepsilon\mathcal{K}(A)$ is homotopically equivalent to the product $_\varepsilon L_0(A) \times B_\varepsilon O(A)^+$ ; a more functorial definition should use the notion of suspension SA of a ring A (see e.g. [1], appendix A).

**1.2.** In a more classical way, the classifying space of algebraic K-theory $\mathcal{K}(A)$ is defined homotopically as the product $K_0(A) \times BGL(A)^+$ and there are two important maps between $\mathcal{K}(A)$ and $_\varepsilon\mathcal{L}(A)$. One is induced by the forgetful functor F from modules with Hermitian forms to modules with no forms ; the other one H from ordinary modules to modules with forms, sending N to H(N), the hyperbolic module associated to N.

$$F : {}_\varepsilon\mathcal{L}(A) \longrightarrow \mathcal{K}(A) \quad \text{and} \quad H : \mathcal{K}(A) \longrightarrow {}_\varepsilon\mathcal{L}(A).$$

Following the notation of [6], we define $_\varepsilon\mathcal{V}(A)$ as the homotopy fiber of F and $_\varepsilon\mathcal{U}(A)$ as the homotopy fiber of H. We write

$$_\varepsilon V_n(A) = \pi_n({}_\varepsilon\mathcal{V}(A)) \text{ and } {}_\varepsilon U_n(A) = \pi_n({}_\varepsilon\mathcal{U}(A)).$$

**1.3. THEOREM** (the fundamental theorem of Hermitian K-theory)[6])**.** *Let us assume that* $1/2 \in A$. *Then*, *there is a natural homotopy equivalence between* $_\varepsilon\mathcal{V}(A)$ *and the loop space of* $_{-\varepsilon}\mathcal{U}(A)$. *In particular,* $_\varepsilon V_n(A) \cong {}_{-\varepsilon}U_{n+1}(A)$. *Moreover, if A is a Banach algebra, the same statement is true for the topological analogs (i.e. replacing* $BGL(A)^+$ *by* $BGL(A)^{top}$ *and* $B_\varepsilon O(A)^+$ *by* $B_\varepsilon O(A)^{top}$).

**1.4.** As it was shown in [4], this theorem, applied to semi-simple finite dimensional algebras over **R**, implies real and complex Bott periodicity [2]. If we neglect 2-torsion phenomenon, theorem 3 is equivalent to the 4-periodicity of the higher "Witt groups"

$$_\varepsilon W_n(A) = \text{Coker } [K_n(A) \longrightarrow {}_\varepsilon L_n(A)].$$

More precisely, let us write **Z'** for the ring **Z**[1/2]. In [6] we proved the existence of elements

$$u_2 \in {}_{-1}W_2(\mathbf{Z'}) \text{ and } u_{-2} \in {}_{-1}W_{-2}(\mathbf{Z'})$$

such that their cup- product is equal to 4 times the unit element of the Witt ring $_1W_0(\mathbf{Z'})$. As a consequence, we have the following theorem.

**1.5. THEOREM.** *Let us assume again that* $1/2 \in A$. *Then there exist* $W_*$-*module maps*

$$\beta : {}_\varepsilon W_n(A) \longrightarrow {}_{-\varepsilon}W_{n+2}(A) \text{ and } \beta' : {}_{-\varepsilon}W_{n+2}(A) \longrightarrow {}_\varepsilon W_n(A)$$

*such that* $\beta.\beta'$ *and* $\beta'.\beta$ *are both multiplication by* 4.

**1.6. Remark.** The Witt groups $W_*(A)$ involved are "$W_*$-modules" in the following sense: they are modules over the ring $_\varepsilon W_*(\mathbf{Z}')$ (all possible $\varepsilon$) via the following cup-product (for any $\mathbf{Z}'$-algebra A):

$$_\eta W_i(A) \times {_\varepsilon W_j(\mathbf{Z}')} \longrightarrow {_{\eta\varepsilon} W_{i+j}(A)}.$$

**1.7.** By considering $(u_2)^4$ and $(u_{-2})^4$, we can define in the same way $W_*$-module maps (in both directions) between for instance $_\varepsilon W_n(A)$ and $_\varepsilon W_{n+8}(A)$, such that their composites are multiplication by $4^4 = 256$, a result which is far from being optimal if we compare it with classical Bott periodicity for Banach algebras. In order to prove the statement quoted in the introduction, we have to choose "better" elements $u_{-\alpha}$ and $u_\alpha$ in $_1W_{-\alpha}(\mathbf{Z}')$ and $_1W_\alpha(\mathbf{Z}')$ respectively (for $\alpha = 8r$ or $8r - 4$).

**1.8. LEMMA.** *There exists an element $u_{-8}$ in $_1W_{-8}(\mathbf{Z}')$ whose image in $_1W_{-8}^{top}(\mathbf{R})$ is $8\,y$, where y is a generator of $_1W_{-8}^{top}(\mathbf{R}) \approx \mathbf{Z}$. In the same way, there exists an element $u_{-4}$ in $_1W_{-4}(\mathbf{Z}')$ whose image in $_1W_{-4}^{top}(\mathbf{R})$ is $2\,z$, where z is a generator of $_1W_{-4}^{top}(\mathbf{R}) \approx \mathbf{Z}$.*

*Proof.* We use the 12 term exact sequence proved in [4] and [6] p. 278 in both the algebraic and topological contexts, where the $k_n$ and the $k'_n$ groups are zero for $n = -1, -3$ and $-5$. We recall that $k_n(A) = H^0(\mathbf{Z}/2\,;\,K_n(A))$ and $k'_n(A) = H^1(\mathbf{Z}/2\,;\,K_n(A))$. If the letter W denotes the higher Witt groups in the algebraic context, we deduce from this exact sequence the following isomorphisms (for the ring $\mathbf{Z}'$)

$$_1W_0 \approx {_{-1}W_{-2}} \approx {_1W_{-4}} \approx {_{-1}W_{-6}} \approx {_1W_{-8}}\,.$$

The situation is quite different for the ring **R** with its usual topology, which we denote by $W^{top}$. In this case we have either isomorphisms or strict inclusions[2] which are multiplication by 2 in the group **Z**:

$$_1W_0 \approx {_{-1}W_{-2}^{top}} \approx {_1W'_{-4}^{top}} \subset {_1W_{-4}^{top}} \approx {_{-1}W'_{-6}^{top}} \subset {_{-1}W_{-6}^{top}} \approx {_1W'_{-8}^{top}} \subset {_1W_{-8}^{top}}$$

It follows that the homomorphisms

$$\mathbf{Z} \approx {_1W_{-8}(\mathbf{Z}')} \longrightarrow {_1W'_{-8}^{top}(\mathbf{R})} \approx \mathbf{Z}$$

$$\mathbf{Z} \approx {_1W_{-4}(\mathbf{Z}')} \longrightarrow {_1W_{-4}^{top}(\mathbf{R})} \approx \mathbf{Z}$$

are multiplication by 8 or 2 respectively. This concludes the proof of lemma 8.

---

[2] We denote by W' the kernel of the rank map.

**1.9. LEMMA.** *For any* $n > 0$, *there exists an element* $u_{4n}$ *in* $_1W_{4n}(\mathbf{Z}')$ *whose image in* $_1W_{4n}^{top}(\mathbf{R}) \approx \mathbf{Z}$ *is twice a generator.*

*Proof.* According to [1], we have a homotopy-cartesian square (where the subscript $_\#$ means 2-adic completion).

$$\begin{array}{ccc} B_1O(\mathbf{Z}')^+_\# & \longrightarrow & BO_\# \times BO_\# \\ \downarrow & & \downarrow \\ B_1O(\mathbf{F}_3)^+_\# & \longrightarrow & BO_\# \end{array}$$

By taking homotopy groups, we find that the map

$$\mathbf{Z} \oplus \mathbf{Z}/2 \approx {}_1W_{4n}(\mathbf{Z}') \approx {}_1W'_{4n}(\mathbf{Z}') \longrightarrow {}_1W'^{top}_{4n}(\mathbf{R}) \approx \mathbf{Z}$$

is an epimorphism. Since the canonical map

$$\mathbf{Z} \approx {}_1W'^{top}_{4n}(\mathbf{R}) \longrightarrow {}_1W^{top}_{4n}(\mathbf{R}) \approx \mathbf{Z}$$

is multiplication by 2, the lemma follows.

**1.10. THEOREM.** *Let* A *be a ring with* 2 *invertible in* A *and* $\alpha$ = 8r *or* 8r - 4 *with* r > 0. *Then, there exist natural homomorphisms of* $W_*$*--modules*

$$\beta : {}_\varepsilon W_i(A) \longrightarrow {}_\varepsilon W_{i+\alpha}(A) \text{ and } \beta' : {}_\varepsilon W_{i+\alpha}(A) \longrightarrow {}_\varepsilon W_i(A)$$

*whose composites are multiplication by* $2.8^r$.

*Proof.* Let us look first at the simple case when 1/2 belongs to A. If $\alpha = 8r$, we define $\beta$ as the cup-product with $u_\alpha$ and $\beta'$ as the cup-product with $(u_{-8})^r$. From lemmas 8 and 9 above, we deduce that their composites are mutiplication by $2.8^r$ as expected. If $\alpha = 8r + 4$, we define again $\beta$ as the cup-product with $u_\alpha$ and we define $\beta'$ as the cup-product with $(u_{-8})^r.(u_{-4})$.
The images of these elements in the topological Witt theory of $\mathbf{R}$ are multiples of the generators of exponents 2 and $2.8^r$ respectively. On the other hand, the cup-product (for any $\alpha \equiv 4$ mod. 8)

$$\mathbf{Z} \times \mathbf{Z} = W^{top}_\alpha(\mathbf{R}) \times W^{top}_{-\alpha}(\mathbf{R}) \longrightarrow W_0(\mathbf{R}) = \mathbf{Z}$$

is well known to be $(x, y) \mapsto 4xy$. Therefore, we finally get the number $4.2.2.8^r = 2.8^{r+1}$ as expected in this case.

**1.11. THEOREM.** *Let* A *be any ring, we may also define homomorphisms*

$$\gamma : {}_\varepsilon W_i(A) \longrightarrow {}_\varepsilon W_{i+\alpha}(A) \text{ and } \gamma' : {}_\varepsilon W_{i+\alpha}(A) \longrightarrow {}_\varepsilon W_i(A)$$

*whose composites are multiplication by* 64 (*resp.* 512) *if* $\alpha = 4$ (*resp.* 8).

*Proof.* Let us consider first any ring B where 2 is a non zero divisor and denote by $\Sigma$ the multiplicative set $(2^n)$. If $K_0(B) = K_{-1}(B) = 0$, there is a localization exact sequence of Witt groups proved in [5] (where W denotes the group ${}_1 W$ and $B_\Sigma$ the localized ring with respect to $\Sigma$)

$$W(B) \longrightarrow W(B_\Sigma) \longrightarrow W(B, \Sigma) \longrightarrow W_{-1}(B) \longrightarrow W_{-1}(B_\Sigma).$$

Here W(B, S) denotes the Witt group of quadratic B-modules of homological dimension $\leq 1$ with S-torsion. We apply this exact sequence to the ring $B = S^\alpha(\mathbf{Z})$, the $\alpha$-suspension of the ring $\mathbf{Z}$, in order to produce the analogs of the elements $u_{-\alpha}$, for $\alpha = 4$ or 8. It is shown in [5] p. 383-385 that the group W(B, S) has at most 8-torsion. It follows that the element $8u_{-\alpha}$ belongs to the image of the map

$$L : W(S^\alpha \mathbf{Z}) \longrightarrow W(S^\alpha \mathbf{Z}').$$

We may now choose an element $v_{-\alpha}$ of $W(S^\alpha \mathbf{Z})$ such that $L(v_{-\alpha}) = 8u_{-\alpha}$.

On the other hand, we have defined in [6], p. 248-251 an element u in ${}_{-1}W_2(\mathbf{Z})$ whose image in ${}_{-1}W_2^{top}(\mathbf{R}) \approx \mathbf{Z}$ is a topological generator. Elementary considerations of topological K-theory [6] show that $(u^2)$ in ${}_1 W_4(\mathbf{Z})$ maps to the double of the generator of ${}_1 W_4^{top}(\mathbf{R}) \approx \mathbf{Z}$. Therefore, if $\alpha = 4$, we see that by taking the cup-products with $v_{-4}$ and $u^2$, we are losing a factor 8 in the previous theorem. In other words, we can define homomorphisms

$$\beta : {}_\varepsilon W_n(A) \longrightarrow {}_\varepsilon W_{n+4}(A) \text{ and } \beta' : {}_\varepsilon W_{n+4}(A) \longrightarrow {}_\varepsilon W_n(A)$$

whose composites are multiplication by $8.8 = 64$.

If $\alpha = 8$, the element $(u)^4$ in ${}_1 W_8(\mathbf{Z})$ maps to 8 y, where y is a generator of the group ${}_1 W_8^{top}(\mathbf{R}) \approx \mathbf{Z}$. Now, if we choose $v_8$ to be $(u)^4$, we see that the image of

$$v_8 \cup v_{-8} \in {}_1 W_0(\mathbf{Z}) = \mathbf{Z} \oplus \text{8-torsion}$$

in ${}_1 W_0(\mathbf{R}) \approx \mathbf{Z}$ is $8.8.8 = 512$, a result which finally proves the last part of the theorem.

## 2. Filtration of the Witt ring W(A) for a commutative ring A.

**2.1.** Let us come back to the simplest case when A is commutative and $1/2 \in A$. Then the images of $_{(-1)^n}W_{2n}(A)$ in $_1W(A) = W(A)$, obtained by the cup-product with the "Bott elements" $(u_{-2})^n = u_{-2n} \in {}_{(-1)^n}W_{-2n}(\mathbf{Z}') = \mathbf{Z}$, form a decreasing filtration $(F_{2n})$ by *ideals* of the Witt *ring* W(A) (the ring structure on W(A) is induced by the ring structure on A). Moreover, we have $F_{2n}F_{2p} \subset F_{2(n+p)}$ as a direct consequence of the definition. On the other hand, we have defined $F_r$ for r odd in the introduction. The fact that

$$F_{2n+2} \subset F_{2n+1} \subset F_{2n}$$

is a consequence of the following exact sequence (with $\varepsilon = (-1)^n$):

(S) $\quad {}_\varepsilon W_{2n+1}(A) \longrightarrow k_{2n+1}(A) \longrightarrow {}_{-\varepsilon}W_{2n+2}(A) \longrightarrow {}_\varepsilon W_{2n}(A)' \longrightarrow k'_{2n+1}(A)$

proved in [6] p. 278. From the definition, it is also clear that $F_{2r+1}.F_{2s} \subset F_{2r+2s+1}$. It is more delicate to prove the inclusion

$$F_{2r+1}.F_{2s+1} \subset F_{2r+2s+2} .$$

Assuming the ring is regular, we can reduce the problem to the case where $r = s = 0$ by considering the loop rings[3] $C = \Omega^r(A)$ and $D = \Omega^s(A)$ as in the "polynomial K or L-theory" defined by Villamayor and the author [4]. We reinterpret $F_{2r+1}$ as the image (in W(A)) of $_\varepsilon W'(C)$ (with $\varepsilon = (-1)^r$) which is the kernel of the "desuspended" rank map :

$$_\varepsilon W(C) \longrightarrow k_0(C) = H^0(\mathbf{Z}/2 ; K_0(C)) = k_{2r}(A).$$

In the same way, we reinterpret $F_{2s+1}$ as the image of $_\eta W'(D)$ (with $\eta = (-1)^s$) which is the kernel of the other desuspended rank map $_\eta W(D) \longrightarrow k_0(D) = H^0(\mathbf{Z}/2 ; K_0(D)) = k_{2s}(A)$. Now we have to show that if $x \in {}_\varepsilon W'(C)$ and $y \in {}_\eta W'(D)$, then the product x.y belongs to the image of the periodicity map (with $\lambda = \varepsilon \eta$):

$$_{-\lambda}W_2(C \otimes_A D) \longrightarrow {}_\lambda W(C \otimes_A D)' = \mathrm{Ker}\,[{}_\lambda W(C \otimes_A D) \longrightarrow k_0(C \otimes_A D)].$$

---

[3] These loop rings (considered as A-algebras) are not unital. One should add a unit and take the associated reduced theory.

Note that for a general ring R, the image of the periodicity map $_{-\lambda}W_2(R) \longrightarrow {_\lambda}W'(R)$ is the kernel of the "discriminant map"

$${_\lambda}W'(R) \longrightarrow H^1(\mathbf{Z}/2 ; K_1(R)) = k'_1(R)$$

as proved using the exact sequence (S) above for n = 0. This discriminant map may be described concretly as follows. Any element x of $_\lambda W'(R)$ is represented by the class of a difference (M, g) - (M, g') where g, g' : M $\longrightarrow$ M* are non degenerate $\lambda$-Hermitian forms. The "discriminant" of x, i.e. its image in $k'_1(R) = H^1(\mathbf{Z}/2 ; K_1(R))$ is the class of the automorphism $g'^{-1}g$. The following lemma follows quite directly from this definition :

**2.2. LEMMA.** *Let* x = (M, g) - (M, g') *and* y = (N, h) - (N, h') *be two elements of* $_\lambda W'(R)$. *Then the product* x.y *in* $_\lambda W'(R)$ *has trivial discriminant. Therefore, x.y belongs to the image of the periodicity map* $_{-\lambda}W_2(R) \longrightarrow {_\lambda}W'(R)$.

*Proof.* We have

$$xy = (M \otimes N, g \otimes h) + (M \otimes N, g' \otimes h') - (M \otimes N, g \otimes h') - (M \otimes N, g' \otimes h)$$

The result then follows from the identity $(g \otimes h)(g' \otimes h')(g \otimes h')^{-1}(g' \otimes h)^{-1} = 1$.

**2.3.** From the previous discussion, il follows that $F_n F_p \subset F_{n+p}$ in general. Conversely, if A is a commutative field of characteristic $\neq 2$ (with trivial involution), it is well known that the elements of $F_2$ belong to $F_1 F_1$. Therefore, if we write as usual I for the ideal $F_1$, we have $I^2 = F_2$
On the other hand, for n = 0, the exact sequene (S) above may be written as follows :

$$_1W_1(A) \longrightarrow k_1(A) \longrightarrow {_{-1}}W_2(A) \longrightarrow {_1}W(A)'$$

Since the first map is surjective, $_{-1}W_2(A)$ may be identified with $I^2$. Therefore, $F_3$ is the subgroup [denoted above by $_{-1}W_2(A)'$] of $_{-1}W_2(A)$ which is the kernel of the rank map
$$_{-1}W_2(A) \longrightarrow k_2(A) = k_2^M(A)$$
Since $k_2^M(A) = I^2/I^3$, we see that we also have $F_3 = I^3$.

**2.4.** For a general commutative ring A as above (again with 2 invertible), we have $I^n \subset F_n$ with a <u>strict</u> inclusion in general. In order to show this last point, we choose the example A = **Q** the field of rational numbers. In this case, according to [3], we have

$$W_8(\mathbf{Q}) = W_8(\mathbf{Z'}) \oplus \text{2-torsion} = \mathbf{Z} \oplus \text{2-torsion}$$

According to [1] again, we know that the image of a free generator of the factor $\mathbf{Z}$ goes to twice a generator of $W_8^{top}(\mathbf{R}) \cong \mathbf{Z}$. Since the Bott element in $W_{-8}(\mathbf{Z'})$ has an image in $W_{-8}^{top}(\mathbf{R})$ which is 8 times a generator, we conclude that $F_8$ is $2^4 \mathbf{Z}$, compared to $I^8 = 2^8 \mathbf{Z}$, where $\mathbf{Z}$ is the free summand of $W(\mathbf{Q}) = \mathbf{Z} \oplus \text{2-torsion}$.

**2.5. THEOREM.** *Let* A *be a commutative ring with* $1/2 \in A$ *and* I *be the augmentation ideal of the Witt ring* W(A)*. Then, if* $m \leq n$*, there is a subquotient of the group* $_{(-1)^m}W_{2m}(A)$ *which is isomorphic to* $I^{2n}$*. In particular, if* $I^{2n} \neq 0$*. we have* $_{(-1)^m}W_{2m}(A) \neq 0$ *for* $m \leq n$.

*Proof.* We have $F_{2m} = \text{Im}\ [\ _{(-1)^m}W_{2m}(A) \longrightarrow W_0(A)\ ]$ which is a quotient of the group $_{(-1)^m}W_{2m}(A)$. Since $I^{2n} \subset F_{2n} \subset F_{2m}$, for all $n \geq m$, there is a subquotient of $_{(-1)^m}W_{2m}(A)$ isomorphic to $I^{2n}$. Therefore, $_{(-1)^m}W_{2m}(A) \neq 0$ if $I^{2n} \neq 0$.

**2.6.** Let us now assume that A is a field with characteristic $\neq 2$ and with the trivial involution. In this case, according to Milnor's conjecture, proved by Voevodsky [7], the quotient group $I^{2n}/I^{2n+1}$ is isomorphic to the Milnor K-group reduced mod. 2 :

$$k_{2n}^M(A) = K_{2n}^M(A)/\ 2\ K_{2n}^M(A)$$

If we assume that this group is $\neq 0$, the same argument shows that $_{(-1)^m}W_{2m}(A) \neq 0$ as above. However, since $k_{2n}^M(A)$ is a 2-group, we can improve this result, at least if $_{(-1)^m}W_{2m}(A) = G$ is finitely generated as an abelian group. In this case, we can write

$$G = (\mathbf{Z})^r \oplus \bigoplus_\alpha (\mathbf{Z}/2^\alpha)^{n_\alpha} \oplus H.$$

where H is a finite group of odd order. With these notations, we shall call $r + \sum_\alpha n_\alpha$ the "2-rank" of G. The following statement is a consequence of Voevodsky's theorem :

**2.7 THEOREM.** *Let* A *be a field of characteristic* $\neq 2$*, provided with the trivial involution. Then*[4]*, for all* $m \leq n$*, there is a subquotient of* $_{(-1)^m}W_{2m}(A)$ *which is isomorphic to* $k_{2n}^M(A)$*. In particular, if* $k_{2n}^M(A)$ *is a* $\mathbf{Z}/2$ *vector space of dimension* r *and if the group* $_{(-1)^m}W_{2m}(A)$ *is finitely generated, the 2-rank of* $_{(-1)^m}W_{2m}(A)$ *is* $\geq r$.

---

[4] The case $m = 0$ is a direct consequence of Milnor's conjecture.

**2.8.** We would like to understand the graded group associated to this new filtration $(F_n)$ of the Witt ring $W(A)$ better. For instance, for the ring $A = \mathbf{Z}' = \mathbf{Z}(1/2)$, we find the following graded groups, starting from $F_0/F_1 = \mathbf{Z}/2$ :

$$\mathbf{Z}/2,\ \mathbf{Z}/2,\ \mathbf{Z}/2,\ \mathbf{Z}/2,\ 0,\ 0,\ 0, 0,\ \mathbf{Z}/2$$

The first 3 groups $\mathbf{Z}/2$ are deduced from the fact that $F_n = I^n$ for $n \leq 3$. In particular $F_3 = 8\mathbf{Z}$. For the other groups, we investigate the maps

$$\mathbf{Z} \cong {}_1W_4(\mathbf{Z}') \longrightarrow {}_1W_0(\mathbf{Z}') \cong \mathbf{Z} \oplus \mathbf{Z}/2 \text{ and } \mathbf{Z} \cong {}_1W_8(\mathbf{Z}') \longrightarrow {}_1W_0(\mathbf{Z}') \cong \mathbf{Z} \oplus \mathbf{Z}/2$$

Since they are both multiplication by 16 on the free summand, we get

$$F_4 = F_5 = F_6 = F_7 = F_8 = 16\,\mathbf{Z}$$

From another viewpoint, the groups $H^n(\mathbf{Z}/2, K_n(\mathbf{Z}')) = H^\mu(\mathbf{Z}/2, K_n(\mathbf{Z}'))$, with $\mu \equiv n \bmod 2$, can be determined from the computations of Rognes and Weibel quoted in [1]. For $n > 0$, the action of $\mathbf{Z}/2$ on $K_n(\mathbf{Z}')$ is $x \mapsto -x$ on the free part and the identity on the 2 primary torsion part. Therefore, we find the following sequence of groups (starting from $n = 1$) :

$$\mathbf{Z}/2,\ \mathbf{Z}/2,\ \mathbf{Z}/2,\ \mathbf{Z}/2, 0, \mathbf{Z}/2, 0, 0,\ \mathbf{Z}/2$$

This sequence is almost the same as the one associated to the filtration $(F_r)$ above. This shows that our filtration (at least in this example) is related to Quillen's K-theory pretty much the same as the filtration by powers of the augmentation ideal I is related to Milnor's K-theory. However, this relation should be more subtle for general rings A. It is probably linked to a (conjectural) spectral sequence ${}_\varepsilon E$ starting with ${}_\varepsilon E^{p,q} = H^p(\mathbf{Z}/2\,;\, K_{-q}(A))$ (with $q \leq 0$) converging to the graded group associated to a certain filtration[5] of ${}_\varepsilon W_n(A)$, with $n = -q - p$. For $A = \mathbf{Z}'$, this spectral sequence should be a consequence of the theorem C proved in [1].

---

[5] which should be the same as the one introduced in this section for $n = 0$.